\input amstex
\documentstyle{amsppt}
\topmatter \magnification=\magstep1 \pagewidth{5.2 in}
\pageheight{6.7 in}
\abovedisplayskip=10pt \belowdisplayskip=10pt
\parskip=8pt
\parindent=5mm
\baselineskip=2pt
\title
 An analogue of Lebesgue-Radon-Nikodym theorem with respect to $p$-adic $q$-invariant
 distribution on $\Bbb Z_p$\endtitle
\author   Taekyun Kim    \endauthor

\affil{ {\it Institute of Science Education,\\
        Kongju National University, Kongju 314-701, S. Korea\\
        e-mail: tkim64$\@$hanmail.net ( or tkim$\@$kongju.ac.kr)}}\endaffil
        \keywords $p$-adic $q$-invariant measure, Lebesgue decomposition,
        $q$-Volkenborn integral, Lebesgue-Radon-Nikodym theorem
Bernoulli numbers
\endkeywords
\thanks  2000 Mathematics Subject Classification:  11S80, 11B68, 11M99 .\endthanks
\abstract{ The purpose of this paper is to derive the analogue of
Lebesgue-Radon-Nikodym  theorem with respect to $p$-adic
$q$-invariant distribution on $\Bbb Z_p$ which is defined by
author in [1].
 }\endabstract
\rightheadtext{   Lebesgue-Radon-Nikodym theorem associated with
$\mu_q$ } \leftheadtext{T. Kim}
\endtopmatter

\document

\head \S 1. Introduction \endhead Let $p$ be a fixed prime.
Throughout this paper $\Bbb Z ,$ $\Bbb Z_p,\,\Bbb Q_p ,$ and $\Bbb
C_p$ will, respectively, denote the ring of rational integers, the
ring of $p$-adic rational integers, the field of $p$-adic rational
numbers and the completion of algebraic closure of $\Bbb Q_p , $
cf.[1, 2, 3].  Let $v_p$ be the normalized exponential valuation
of $\Bbb C_p$ with $|p|=p^{-v_p(p)}=p^{-1}$ and let $a+p^N\Bbb
Z_p=\{x\in\Bbb Z_p| x \equiv a (\mod p^N)\},$ where $a\in\Bbb Z$
lies in $0\leq a<p^N .$ In this paper we assume that $q\in\Bbb
C_p$ with $|1-q|<p^{-\frac{1}{p-1}}$ as an indeterminate. We now
use the notation
$$[x]_q=[x:q]=\frac{1-q^x}{1-q}, \text{ cf. [1, 2, 3]}.$$
For any positive integer $N$ we set
$$\mu_q(a+p^N\Bbb Z_p=\frac{q^a}{[p^N]_q }, \text{ ( see [1] )},\tag 1
$$ and this can be extended to a distribution on $\Bbb Z_p$. In
this paper $\mu_q$ will be called $p$-adic $q$-invariant
distribution on $\Bbb Z_p $.

 For $f\in C^{(1)}(\Bbb Z_p , \Bbb
C_p)=\{f|f:\Bbb Z_p \rightarrow \Bbb C_p \text{ is
$C^{(1)}$-function } \},$ the above distribution $\mu_q$ yields an
integral  as follows:
$$I_q(f)=\int_{\Bbb Z_p}f(x) d\mu_q(x)=\lim_{N\rightarrow \infty}
\frac{1}{[p^N]_q}\sum_{x=0}^{p^N-1}f(x)q^x , \text{ ( see [1] ) .
}\tag 2$$ This integral is called $q$-Volkenborn integration,
cf.[1]. In [1, 3], the $q$-Bernoulli numbers are defined as
$$\beta_{m,q}=\int_{\Bbb Z_p}q^{-x}[x]_q^m d\mu_q(x),  \text{
$m\geq 0$ }. \tag 3$$ Note that
$$\beta_{0,q}=\frac{q-1}{\log q}, \text{ }
\beta_{m,q}=\frac{1}{(1-q)^m}\sum_{i=0}^{\infty}\binom mi
(-1)^{m-i}\frac{i}{[i]_q}, \text{ cf.[1]}. $$ Also, the
$q$-Bernoulli polynomials are defined by
$$\beta_{n,q}(x)=\int_{\Bbb Z_p}q^{-y}[x+y]_q^n
d\mu_q(y)=\sum_{k=0}^n\binom nk q^{kx}\beta_{k,q}[x]_q^{n-k},
\text{ $n\geq 0$}. \tag 4$$ Let $C(\Bbb Z_p , \Bbb C_p)$ be the
space of continuous function on $\Bbb Z_p $ with values in $\Bbb
C_p ,$ provided with norm  $||f||_{\infty}=\sup_{x\in\Bbb Z_p}
|f(x)|. $ The difference quotient $\Delta_1 f$ of $f$ is the
function of two variables given by
$$\Delta_1 f(m,x)=\frac{f(x+m)-f(x)}{m}, \text{ for all $x,
m\in\Bbb Z_p$, $m\neq 0$.}$$ A function $f:\Bbb Z_p \rightarrow
\Bbb C_p $ is said to be a Lipschitz function if there exists  a
constant $M >0$ ( the Lipschitz constant  of $f$ ) such that
$$|\Delta_1 f(m,x)|\leq M \text{ for all $m\in\Bbb Z_p\setminus
\{0\} $ and $x\in\Bbb Z_p$.}$$ The $\Bbb C_p$-linear space
consisting of all Lipschitz function ( or $C^{(1)}$-function ) is
denoted by $Lip(\Bbb Z_p, \Bbb C_p)$ ( or $C^{(1)}(\Bbb Z_p, \Bbb
C_p $)). This space is a Banach space with respect to the norm
$||f||_1 =||f||_{\infty}\vee ||\Delta_1 f||_{\infty}. $ The
purpose of this paper is to derive  the analogue of
Lebesgue-Radon-Nikodym theorem  with respect to $p$-adic
$q$-invariant distribution on $\Bbb Z_p$ which is defined by
author in [1].

 \head 2. Analogue of Lebesgue-Radon-Nikodym theorem with respect to $\mu_q$ on $\Bbb Z_p$  \endhead

For $f\in C^{(1)}(\Bbb Z_p, \Bbb C_p)$, the $q$-Volkenborn
integral is defined by
$$ I_q(f)=\int_{\Bbb Z_p} f(x) d\mu_q(x)=\lim_{N\rightarrow
\infty}\frac{1}{[p^N]_q}\sum_{x=0}^{p^N-1} f(x)q^x, \text{ cf.
[1].}$$ This definition is based on the $p$-adic $q$-distribution
on $\Bbb Z_p$, defined by $\mu_q(a+p^N\Bbb
Z_p)=\frac{q^a}{[p^N]_q}.$ By the meaning of the extension of
$q$-Volkenborn integral, we consider the below weakly (strongly)
$p$-adic $q$-invariant distribution $\mu_q$ on $\Bbb Z_p$,
satisfying
$$\left|[p^n]_q\mu_q(a+p^n\Bbb Z_p)-[p^{n+1}]_q\mu_q(a+p^{n+1}\Bbb
Z_p) \right|\leq \delta_n, \tag 5$$ where $\delta_n\rightarrow 0$,
$a\in\Bbb Z$, and $\delta_n$ is independent of $a$ ( for strongly
$p$-adic $q$-distribution, $\delta_n$ is replaced by $cp^{-n}$,
where $c$ is positive real constant ). Let $f(x)$ be be a function
on $\Bbb Z_p$. Then the $q$-Volkenborn integral of $f$ with
respect to weakly $\mu_q$ is
$$\int_{\Bbb Z_p}f(x) d\mu_q(x)=\lim_{N\rightarrow \infty}
\sum_{x=0}^{p^N-1}f(x)\mu_q(x+p^N\Bbb Z_p), \tag 6$$ if the limit
exists. In this case we say that $f$ is $q$-Volkenborn integrable
with respect to $\mu_q$. A distribution $\mu_q$ on $\Bbb Z_p$ is
called $1$-admissible if $$\left| \mu_q(a+p^N\Bbb Z_p) \right|\leq
\frac{1}{|[p^N]_q|}c_n, $$ where $c_n\rightarrow 0.$ If $\mu_q$ is
1-admissible on $\Bbb Z_p$, then we see that $\mu_q$ is a weakly
$p$-adic $q$-invariant distribution on $\Bbb Z_p$. Let $\mu_q$ be
a weakly $p$-adic $q$-invariant distribution on $\Bbb Z_p$, then
$\lim_{N}[p^N]_q\mu_q(x+p^N\Bbb Z_p)$ exists for all $x\in\Bbb
Z_p, $ and converges uniformly with respect to $x$. We define the
$q$-analogue of Radon-Nikodym derivative of $\mu_q$ with respect
to $p$-adic $q$-invariant distribution on $\Bbb Z_p$ as follows:
$$f_{\mu_q}(x)=\lim_{N\rightarrow \infty}[p^N]_q\mu_q(x+p^N\Bbb
Z_p). \tag 7$$ Let $x, y \in\Bbb Z_p$ with $|x-y|_p=p^{-m}.$ Then
$x+p^m\Bbb Z_p=y+p^m\Bbb Z_p.$ Hence, we note that
$$\aligned
&|f_{\mu_q}(x)-f_{\mu_q}(y)|=|f_{\mu_q}(x)-[p^n]_q\mu_q(x+p^n\Bbb
Z_p)-f_{\mu_q}(y)+[p^n]_q\mu_q(y+p^n\Bbb Z_p) +\\
&[p^n]_q\mu_q(x+p^n\Bbb Z_p)-[p^m]_q\mu_q(x+p^m\Bbb
Z_p)-[p^n]_q\mu_q(y+p^n\Bbb Z_p)+[p^m]_q\mu_q(y+p^m\Bbb Z_p)|.
\endaligned \tag 8$$
If $\mu_q$ is strongly $p$-adic $q$-invariant distribution on
$\Bbb Z_p$, then we see that
$$|[p^n]_q\mu_q(x+p^n\Bbb Z_p)-[p^m]_q\mu_q(x+p^m\Bbb Z_p)|\leq
Cp^{-m}, \tag 9$$ for some constant $C$, and $m<n .$

For $n >> 0$, we note that
$$|f_{\mu_q}(x)-f_{\mu_q}(y)|\leq C_1 p^{-m}=C_1|x-y|,  $$
where $C_1$ is positive real constant.

Therefore we obtain the following proposition :
\proclaim{Proposition 1} Let $\mu_q$ be a strongly $p$-adic
$q$-invariant distribution on $\Bbb Z_p$. Then $f_{\mu_q}\in
Lip(\Bbb Z_p, \Bbb C_p).$
\endproclaim
Let $f\in C^{(1)}(\Bbb Z_p, \Bbb C_p)$. For any positive integer
$a$, with $ a< p^N,$ define
$$\mu_{f,q}(a+p^n\Bbb Z_p)=\int_{a+p^n\Bbb Z_p}q^{-xp^n-a}f(x)
d\mu_q(x), \tag 10 $$ where the integral is $q$-Volkenborn
integral.

Note that
$$\aligned
\mu_{f,q}(a+p^n \Bbb Z_p)&=\lim_{m\rightarrow
\infty}\frac{1}{[p^{m+n}]_q}\sum_{x=0}^{p^m-1}f(a+p^n
x)=\lim_{m\rightarrow
\infty}\frac{1}{[p^m]_q}\sum_{x=0}^{p^{m-n}-1}f(a+p^n
x)\\
&=\frac{1}{[p^n]_q}\int_{\Bbb Z_p} q^{-p^n x}f(a+p^n x)
d\mu_{q^{p^n}}(x).
\endaligned$$
By(10), we obtain the following proposition:

\proclaim{ Proposition 2} For $f, g\in C^{(1)}(\Bbb Z_p, \Bbb C_p
), $ we have
$$\aligned
&\mu_{\alpha f +\beta g ,q}=\alpha \mu_{f,q}+\beta\mu_{g,q}, \\
&|\mu_{f,q}(a+p^n \Bbb Z_p)|\leq ||f||_1 |\frac{1}{[p^n]_q}|.
\endaligned \tag11$$
\endproclaim

Let $P(x)\in\Bbb C_p[[x]_q]$ be an arbitrary polynomials. We now
give the proof that $\mu_{P,q}$ is a strongly $p$-adic
$q$-invariant distribution on $\Bbb Z_p$. It is enough to prove
the statement for $P(x)=[x]_q^k.$ Let $a$ be an integer with
$0\leq a <p^n .$ Then we see that
$$\aligned
&\mu_{P,q}(a+p^n\Bbb
Z_p)=\lim_{m}\frac{1}{[p^m]_q}\sum_{i=0}^{p^{m-n}-1}[a+ip^n]_q^k\\
&=\lim_{m\rightarrow \infty}\frac{1}{[p^m]_q}(
p^{m-n}[a]_q^k+k[a]_q^{k-1}q^a[p^n]_q\sum_{i=0}^{p^{m-n}-1}[i]_{q^{p^n}}+\cdots
+[p^n]_q^kq^{ak}\sum_{i=0}^{p^{m-n}-1}[i]_{q^{p^n}}^k
).\endaligned$$ Because
$[a+ip^n]_q=([a]_q+q^a[p^n]_q[i]_{q^{p^n}}).$
 From (3), we note that
 $$\beta_{n,q}=\int_{\Bbb Z_p}q^{-x}[x]_q^n
 d\mu_q(x)=\lim_{m\rightarrow\infty}\frac{1}{[p^m]_q}\sum_{i=0}^{p^{m}-1}[i]^n.$$
Thus, we have
$$\aligned
&[p^n]_q\mu_{P,q}(a+p^n\Bbb Z_p)\\
&=-\frac{1-q^{p^n}}{p^n\log q}[a]^k+[p^n]_qq^a\binom k1
[a]_q^{k-1}\beta_{1,q^{p^n}}+\cdots+[p^n]_q^kq^{ak}\beta_{k,q^{p^n}}\\
&=\sum_{i=0}^{\infty}\frac{P^{(i)}(a)}{i!}\beta_{i,q^{p^n}}[p^n]_q^i
q^{ai}, \text{ where $P^{(i)}(a)=\left(\frac{d}{d[x]_q}\right)^i
P(x) |_{x=a}.$ }
\endaligned \tag12$$

The Eq.(12) can be rewritten  as
$$\mu_{P,q}(a +p^n\Bbb Z_p)\equiv
\frac{P(a)}{[p^n]_q}\beta_{0,q^{p^n}}+q^a\beta_{1,q^{p^n}}P^{\prime}(a)
\text{ ( $\mod [p^n]_q $) }. \tag13$$

Let $x$ be an arbitrary in $\Bbb Z_p$, $x\equiv x_n$ ($\mod p^n$),
$x\equiv x_{n+1}$ ($\mod p^{n+1}$), where $x_n, x_{n+1}$ are
positive integers such that $0\leq x_n <p^n$ and $0\leq
x_{n+1}<p^{n+1} .$ Then
$$\aligned
&|[p^n]_q\mu_{P,q}(x+p^n\Bbb
Z_p)-[p^{n+1}]_q\mu_{P,q}(x+p^{n+1}\Bbb Z_p)|\\
&=|\sum_{i=0}^{\infty}\frac{P^{(i)}(x_n)}{i!}\beta_{i,q^{p^n}}[p^n]_q^iq^{ai}-\sum_{i=0}^{\infty}\frac{P^{(i)}(x_{n+1})}{i!}
\beta_{i,(q^p)^{p^n}}[p^n]_{q^p}^i[p]_q^i q^{ai}|.
\endaligned\tag14$$
By (3), (8), (13) and (14), we easily see that
$$|[p^n]_q\mu_{P,q}(x+p^n\Bbb
Z_p)-[p^{n+1}]_q\mu_{P,q}(x+p^{n+1}\Bbb Z_p)|\leq C_2p^{-n},
\text{ for some constant $C_2$. } $$ Thus, we note that
$$\aligned
&f_{\mu_{P,q}}(a)=\lim_{n}[p^n]_q\mu_{P,q}(a+p^n\Bbb Z_p)\\
&=\lim_{n\rightarrow \infty}(\frac{q^{p^n}-1}{\log
q^{p^n}}[a]_q^k+[p^n]_qq^a[a]_q^{k-1}\binom
k1\beta_{1,q^{p^n}}+\cdots+[p^n]_q^kq^{ka}\beta_{k,q^{p^n}
})=[a]_q^k=P(a),\\
\endaligned$$
because $\lim_{n\rightarrow \infty}q^{p^n}=1$ for
$|1-q|_p<p^{-\frac{1}{p-1}}.$

For all $x\in\Bbb Z_p ,$ we have $f_{\mu_{P,q}}(x)=P(x)$, since
$f_{\mu_{P,q}}(x)$ is continuous in (7). Now let $g\in
C^{(1)}(\Bbb Z_p, \Bbb C_p) .$ By (12), we easily see that
$$\aligned
\int_{\Bbb Z_p}
&g(x)d\mu_{P,q}(x)\\
&=\lim_{n}\sum_{i=0}^{p^n-1}g(i)\mu_{P,q}(i+p^n\Bbb
Z_p)=\lim_{n}\frac{1}{[p^n]_q}\sum_{i=0}^{p^n-1}g(i)\sum_{j=0}^k\binom
kj \beta_{j,q^{p^n}}[p^n]_q^jq^{ij}[i]_q^{k-j}
\\
&=\sum_{j=0}^k\binom
kj\left(\lim_{n}\beta_{j,q^{p^n}}\frac{[p^n]_q^j}{[p^n]_q}
\sum_{i=0}^{p^n-1}g(i)q^{ij}[i]_q^{k-j}\right)=\lim_{n}\beta_{0,q^{p^n}}
\frac{1}{[p^n]_q}\sum_{i=0}^{p^n-1}g(i)[i]_q^k\\
&=\int_{\Bbb Z_p}g(x)[x]_q^kq^{-x}d\mu_q(x),
\endaligned$$
where the last integral is  $q$-Volkenborn integral on $\Bbb Z_p$
which is defined by author in [1]. Therefore we obtain the
following theorem:
 \proclaim{Theorem 3}
Let $P(x)\in\Bbb C_p [[x]_q]$ be an arbitrary polynomials. Then
$\mu_{P,q}$ is a strongly $p$-adic $q$-invariant distribution, and
for all $x\in\Bbb Z_p,$
$$f_{\mu_{P,q}}(x)=P(x). \tag15$$
Furthermore, for any $g\in C^{(1)}(\Bbb Z_p, \Bbb C_p ),$
$$\int_{\Bbb Z_p}g(x) d\mu_{P,q}(x)=\int_{\Bbb Z_p}g(x)P(x)q^{-x}
d\mu_q(x), \tag16 $$ where the second integral is $q$-Volkenborn
integral.
\endproclaim
Let $f(x)=\sum_{n=0}^{\infty}a_{n,q}\binom xn_q$ be the $q$-Mahler
expansion of the $C^{(1)}$-function of $f$. Here, $\binom
xn_q=\frac{[x]_q[x-1]_q\cdots[x-n+1]_q}{[n]_q[n-1]_q\cdots[2]_q[1]_q},$
cf.[2]. Then we note that $\lim_{n}n|a_{n,q}|=0, $ (see [2]).
Remark. In the recent, the $q$-Mahler expansion was introduced by
K. Conrad, cf.[2].

Consider $f_m(x)=\sum_{i=0}^{m}a_{i,q}\binom xi_q\in\Bbb
C_p[[x]_q].$ Then $||f-f_m||_1\leq \sup_{n\geq m} n|a_{n,q}|, $
cf.[1, 2]. Writing $f=f_m+(f-f_m)$. By using (11), we can observe
that
$$\aligned
&|[p^n]_q\mu_{f,q}(a+p^n\Bbb
Z_p)-[p^{n+1}]_q\mu_{f,q}(a+p^{n+1}\Bbb Z_p)|\\
&\leq \max\{ |[p^n]_q\mu_{f_m,q}(a+p^n\Bbb Z_p) -[p^{n+1}]_q
\mu_{f_m,q}(a+p^{n+1}\Bbb Z_p)|,\\
& |[p^n]_q\mu_{f-f_m , q}(a+p^n\Bbb Z_p)-[p^{n+1}]_q\mu_{f-f_m ,
q}(a+p^{n+1}\Bbb Z_p)| \}.
\endaligned$$
From Theorem 3, we note that $\mu_{f_m,q}$ is a strongly $p$-adic
$q$-invariant distribution on $\Bbb Z_p .$ Moreover, if $m >>0$,
then we know that
$$|[p^n]_q\mu_{f-f_m,q}(a+p^n\Bbb Z_p)|\leq ||f-f_m||_1 \leq
C_3p^{-n} ,$$ where $C_3$ is some positive real constant. Also for
$m>>0$, it follows that $||f||_1=||f_m||_1$ and so
$$|[p^n]_q\mu_{f_m,q}(a+p^n\Bbb Z_p)-[p^{n+1}]_q\mu_{f_m,
q}(a+p^{n+1}\Bbb Z_p)|\leq C_4p^{-n}, $$ where $C_4$ is also some
positive real constant.

We now observe that
$$\aligned
&|f(a)-[p^n]_q\mu_{f,q}(a+p^n\Bbb Z_p)|\\
&=|f(a)-f_m(a)+f_m(a)-[p^m]_q\mu_{f_m, q}(a+p^n\Bbb
Z_p)-[p^n]_q\mu_{f-f_m, q}(a+p^n\Bbb Z_p)|\\
&\leq \max\{|f(a)-f_m(a)|, |f_m(a)-[p^n]_q\mu_{f_m,q}(a+p^n\Bbb
Z_p)|,|[p^n]_q\mu_{f-f_m, q}(a+p^n\Bbb Z_p)|\}\\
&\leq\max\{|f(a)-f_m(a)|, |f_m(a)-[p^n]_q\mu_{f_m,q}(a+p^n\Bbb
Z_p)|,||f-f_m||_1\}.
\endaligned$$
Now if we fix $\epsilon >0 ,$ and fix $m$ such that $||f-f_m||_1
\leq \epsilon, $ then for $ n>>0$, we have
$$|f(a)-[p^n]_q\mu_{f,q}(a+p^n\Bbb Z_p)|\leq \epsilon .$$
That is,
$$f_{\mu_{f,q}}(a)=\lim_{n\rightarrow
\infty}[p^n]_q\mu_{f,q}(a+p^n\Bbb Z_p)=f(a). \tag17$$ Let $m$ be
the sufficiently large number such that
 $||f-f_m||_1 \leq p^{-2n}.$ Then we see that
 $$|\mu_{f-f_m, q}(a+p^n\Bbb Z_p)|\leq \frac{1}{|[p^n]_q|}||f-f_m||_1\leq[p^n]_qp^{-2n}\leq p^{-n}.$$
Hence, we have
$$\aligned
&\mu_{f,q}(a+p^n\Bbb Z_p)=\mu_{f_m,q}(a+p^n\Bbb
Z_p)+\mu_{f-f_m,q}(a+p^n\Bbb Z_p)\equiv \mu_{f_m,q}(a+p^n\Bbb
Z_p)\\
&\equiv
\frac{f_m(a)}{[p^n]_q}\beta_{0,q^{p^n}}-\beta_{1,q^{p^n}}f_m^{\prime}(a)
=\frac{f(a)}{[p^n]_q}\beta_{0,q^{p^n}}-\beta_{1,q^{p^n}}f^{\prime}(a)
\text{ ($mod[p^n]_q$) }.
\endaligned\tag18$$
Let $g\in C^{(1)}(\Bbb Z_p, \Bbb C_p).$ Then
$$\aligned
\int_{\Bbb
Z_p}g(x)d\mu_{f,q}(x)&=\lim_{n}\sum_{i=0}^{p^n-1}g(i)\mu_{f,q}(i+p^n\Bbb
Z_p)\\
&=\lim_{n}\sum_{i=0}^{p^n-1}g(i)\left(\frac{f(i)}{[p^n]_q}\beta_{0,q^{p^n}}-\beta_{1,q^{p^n}}f^{\prime}(i)\right)\\
&=\lim_{n}\left(\frac{1}{[p^n]_q}\sum_{i=0}^{p^n-1}g(i)f(i)-\frac{1}{2}\sum_{i=0}^{p^n-1}g(i)f^{\prime}(i)\right)\\
&=\int_{\Bbb Z_p}g(x)f(x)q^{-x}d\mu_q(x),
\endaligned$$
because $\lim_{n\rightarrow \infty}\beta_{k,q^n}=B_k$, where $B_k$
are the ordinary $k$-th Bernoulli numbers.

Let $\delta$ be the function from $C^{(1)}(\Bbb Z_p, \Bbb C_p)$ to
$Lip(\Bbb Z_p, \Bbb C_p).$ We know that $\mu_q$ is a strongly
$p$-adic $q$-invariant distribution on $\Bbb Z_p$, then
$|f_{\mu_q}(a)-[p^n]_q\mu_q(a+p^n\Bbb Z_p)|\leq Cp^{-n}$ for any
positive integer $n$. Let  $\mu_{1,q}$ be the strongly $p$-adic
$q$-invariant distribution associated to $f_{\mu_q}$(
$\mu_{1,q}=\delta(f_{\mu_q})$). If $\mu_{1,q}$ is the associated
strongly $p$-adic $q$-invariant distribution on $\Bbb Z_p$, then ,
by (18), we have
$$ |[p^n]_q\mu_{1,q}(a+p^n\Bbb Z_p)-f_{\mu_q}(a)|\leq
\max\{p^{-n}, p^{-n}|\beta_{0,q^{p^n}}f^{\prime}(a)|\}\leq
C_5p^{-n}, \text{ for $n>>0$, } $$ where $C_5$ is a some positive
real constant. Because $[p^n]_q\equiv 0$ ($\mod p^n$).

This shows that for sufficient large $n$
$$\aligned
&|\mu_q(a+p^n\Bbb Z_p)-\mu_{1,q}(a+p^n\Bbb Z_p)|\\
&=\frac{1}{|[p^n]_q|}|[p^n]_q\mu_q(a+p^n\Bbb
Z_p)-f_{\mu_q}(a)+f_{\mu_q}(a)-[p^n]_q\mu_{1,q}(a+p^n\Bbb Z_p)|\\
&\leq \max \{p^n p^{-n} C_6, p^np^{-n} C_7\}\leq M,
\endaligned$$
where $C_6, C_7, M$ are some positive real constants.

Thus, $\mu_q-\mu_{1,q}$ is a distribution and bounded. Hence,
$\mu_q-\mu_{1,q}$ is a measure on $\Bbb Z_p.$ Therefore we obtain
the following theorem: \proclaim{ Theorem 4} Let $\mu_q$ be a
strongly $p$-adic $q$-invariant distribution on $\Bbb Z_p$, and
assume that the Radon-Nikodym derivative  $f_{\mu_q}$ on $\Bbb
Z_p$ is a $C^{(1)}$-function. Suppose that $\mu_{1,q}$ is the
strongly $p$-adic $q$-invariant distribution associated to
$f_{\mu_q}$ $( \mu_{1,q}=\delta (f_{\mu_q}))$, then there exists a
measure $\mu_{2,q}$ on $\Bbb Z_p$ such that
$$\mu_q=\mu_{1,q}+\mu_{2,q}. \tag 19$$
\endproclaim

Eq.(19) is $p$-adic analogue of Lebesgue decomposition with
respect to strongly  $p$-adic $q$-invariant distribution on $\Bbb
Z_p$.

\enddocument